\documentclass{amsart}

\usepackage{amssymb}

\newtheorem{thm}{Theorem}%[section]

\newtheorem{lem}[thm]{Lemma}
\newtheorem{cor}[thm]{Corollary}

\newtheorem{prop}[thm]{Proposition}

\theoremstyle{definition}
\newtheorem{defn}[thm]{Definition}

\newtheorem{say}[thm]{}
\newtheorem{exmp}[thm]{Example}

\newtheorem{rem}[thm]{Remark}          
\newtheorem{ack}{Acknowledgments}        
\newtheorem{notation}[thm]{Notation}   
  
\newtheorem{defn-thm}[thm]{Definition--Theorem}  %!!!!!!!!!!!!!!!!!!!!!!!!
\theoremstyle{remark}
\newtheorem{complem}[thm]{Complement}%!!!!!!!!!!!!!!!!!!!!!!
\newtheorem{claim}[thm]{Claim}

\renewcommand{\c}[0]{{\mathbb C}}  

\renewcommand{\o}[0]{{\mathcal O}} 
\newcommand{\z}[0]{{\mathbb Z}}

\newcommand{\p}[0]{{\mathbb P}}

\newcommand{\q}[0]{{\mathbb Q}}

\newcommand{\qtq}[1]{\quad\mbox{#1}\quad}

\newcommand{\rank}[0]{\operatorname{rank}}
\newcommand{\mult}[0]{\operatorname{mult}}

\newcommand{\lcm}[0]{\operatorname{lcm}}

\def\into{\DOTSB\lhook\joinrel\to}

\begin{document}
\bibliographystyle{amsalpha}

\title{Einstein metrics on connected sums of $S^2\times S^3$}
\author{J\'anos Koll\'ar}
\address{Princeton University,
Princeton, NJ 08544-1000.}
\email{kollar@math.princeton.edu}

%\today
\maketitle

It is still very poorly understood which 5-manifolds carry
an Einstein metric with positive constant.
By Myers' theorem, the fundamental group of such a manifold 
 is finite, 
therefore it is reasonable to concentrate on
the simply connected case.
The most familiar examples are 
 connected sums
of $k$ copies of $S^2\times S^3$.

For $k\leq 9$, 
Einstein metrics on these were constructed
by Boyer, Galicki and Nakamaye \cite{bgn02, bgn03, bg03}.  
In this paper we extend their result to any $k$.

\begin{thm} For every $k\geq 6$, there are 
 infinitely many  $(2k-2)$-dimensional families of 
Einstein metrics on the connected sum
of $k$ copies of $S^2\times S^3$.
\end{thm}

It was known earlier that these spaces carry metrics of
positive Ricci curvature; this is a special case of the
results of \cite{sha-ya}. Sasakian metrics of
positive Ricci curvature on $k\#(S^2\times S^3)$
are constructed in \cite{bgn03c}.
\medskip

The constructions in \cite{bgn02, bgn03, bg03}
exhibit suitable links of singular hypersurfaces 
$0\in Y\subset \c^m$ with $\c^*$-action.
These links are Seifert bundles over the corresponding
weighted projective hypersurfaces $(Y\setminus\{0\})/\c^*$. 
The Einstein metrics
are then constructed from  K\"ahler--Einstein metrics
on these weighted projective hypersurfaces.

Here we look at this construction from the other end.
Starting with a projective variety $X$, we study 
Seifert bundles $Y\to X$. For $X$ smooth, these are described
in \cite{or-wa}. When $X$ is a surface, the topology of $Y$
can be understood   well enough. 
The freedom we gain is that one can start with an arbitrary
algebraic surface, not just with a hypersurface
in a weighted projective space.

The conditions for the method to work are somewhat delicate, but
  there are probably many more
examples obtained similarly.

The natural setting of this approach is to start with an
 orbifold $X$. The main ideas are similar, but this requires
a somewhat lengthy
study of Seifert bundles over orbifolds, which will be done elsewhere.

\section{Seifert bundles over complex manifolds}

Seifert bundles over complex manifolds were introduced and studied in
\cite{or-wa}. The quickest approach is to define them first locally
by explicit formulas and then to patch the local forms together.
(See \cite{or-wa} for a conceptually better, though equivalent, definition.)

\begin{defn} \label{stand.seif.defn}
Let $D^n_{\mathbf z}\subset \c^n$ be the unit polydisc
with coordinates $z_1,\dots,z_n$. Let 
$G_t$ be either $\c^*$ or the unit circle
$S^1\subset \c$, both  with
coordinate $t$. Pick pairwise relatively prime natural numbers $a_1,\dots,a_n$
and
for every $i$ pick  $0<b_i\leq a_i$ such that
$b_i$ is relatively prime to $a_i$. Set $a=\prod a_i$.

Let $\epsilon_i$ be a primitive $a_i$th root of unity
and consider the  $\z/a$ action on $D^n_{\mathbf z}$ given by
$$
\phi: (z_1,\dots,z_n)\mapsto 
(\epsilon_1z_1,\dots,\epsilon_nz_n),
$$
and its lifting to $G_t\times D^n_{\mathbf z}$
$$
\Phi: (t,z_1,\dots,z_n)\mapsto 
((\textstyle{\prod} \epsilon_i^{b_i})t, \epsilon_1z_1,\dots,\epsilon_nz_n).
$$
It is easy to see that 
$$
D^n_{\mathbf z}/\langle\phi\rangle\cong D^n_{\mathbf x},
\qtq{where}
x_i=z_i^{a_i}.
$$
The quotient map $D^n_{\mathbf z}\to D^n_{\mathbf x}$
ramifies along $(z_i=0)$ with multiplicity $a_i$.
Furthermore, $\prod \epsilon_i^{b_i}$ is a primitive  $a$th root of unity,
which implies that $G_t\times D^n_{\mathbf z}/\langle\Phi\rangle$
 is a smooth manifold.
The second projection of $G_t\times D^n_{\mathbf z}$ descends to a map
$$
f:G_t\times D^n_{\mathbf z}/\langle\Phi\rangle\to
 D^n_{\mathbf z}/\langle\phi\rangle\cong D^n_{\mathbf x}.
$$
This is called a {\it standard 
Seifert $G$-bundle} over $D^n_{\mathbf x}$ with  orbit invariants
$(a_i,b_i)$ along $(x_i=0)$.
($(\alpha_i,\beta_i)$ in the notation of \cite{or-wa}.)

Notice that the bounded and $\z/a$-invariant  holomorphic sections of
$\c_t\times D^n_{\mathbf z}\to D^n_{\mathbf z}$ are of the form
$$
t=(\textstyle{\prod} z_i^{b_i})\cdot h(z_1^{a_1},\dots,z_n^{a_n}),
\qtq{where $h$ is holomorphic.}
$$
Thus the  bounded holomorphic sections of
$\c_t\times D^n_{\mathbf z}/\langle\Phi\rangle\to D^n_{\mathbf x}$
form a locally free sheaf whose generator can be thought of as
$\prod x_i^{b_i/a_i}$.
\end{defn}

\begin{defn}\label{seif.S1.defn}
 Let $X$ be a complex manifold and
$D_i\subset X$ smooth divisors intersecting transversally.
For every $i$ pick natural numbers $0<b_i\leq a_i$ such that
$a_i$ and $b_i$ are relatively prime for every $i$.
Assume that $a_i$ and $a_j$ are relatively prime whenever 
$D_i\cap D_j\neq \emptyset$.

Every $x\in X$ has a neighborhood
$x\in U\subset X$ 
and a biholomorphism $\tau:U\cong D^n$ such that every $D_i\cap U$ is 
mapped to a
coordinate hyperplane in $D^n$ by $\tau$ (or the intersection is empty).
This assigns numbers $(a_i,b_i)$ to every coordinate on $D^n$.
(We set $a_j=b_j=1$ for those coordinates that do not correspond to a $D_i$.)

A  {\it Seifert $G$-bundle}
 over $X$ with {\it orbit invariants}
$(a_i,b_i)$ along $D_i$ is a real manifold $L$ with a differentiable
$G$-action  and a 
differentiable  map $f:L\to X$ such that for every neighborhood $U$ 
as above,
 $\tau\circ f:f^{-1}(U)\to U\cong D^n$
is fiber preserving equivariantly diffeomorphic to the corresponding 
standard Seifert model.

For any $P\in X$, the number $\prod_{P\in D_i}a_i$ is
called the {\it multiplicity} of the Seifert fiber over $P$.

A   Seifert $S^1$-bundle is also called  a {\it Seifert bundle},

There is a one--to--one correspondence between
 Seifert $S^1$-bundles and  Seifert $\c^*$-bundles over $X$.
\end{defn}

\begin{defn}\label{seif.C*.defn}
Analogously, a  {\it holomorphic Seifert $\c^*$-bundle} 
over $X$ with  orbit invariants
$(a_i,b_i)$ along $D_i$ is a complex manifold $Y$ with a holomorphic
$\c^*$-action  and a holomorphic map $f:Y\to X$ 
such that for every neighborhood $U$ 
as above,
 $\tau\circ f:f^{-1}(U)\to U\cong D^n$
is fiber preserving equivariantly biholomorphic to the corresponding 
standard Seifert model.

From (\ref{stand.seif.defn}) we see that
$$
X\supset V\mapsto 
\{\mbox{bounded holomorphic sections of $f$ over $V\setminus\cup D_i$}\}
$$
defines a locally free sheaf, denoted by   $B_Y$.
\end{defn}

\begin{say}[Construction of Seifert bundles]\cite[3.9]{or-wa}
\label{seif.contr}
 Let $X$ be a complex manifold 
such that $H_1(X,\z)=0$ .
Assume that we are given
\begin{enumerate}
\item smooth divisors
$D_i\subset X$  intersecting transversally,
\item natural numbers
 $0<b_i\leq a_i$   such that 
\begin{enumerate}
\item  $a_i$ and $b_i$ are relatively prime for every $i$, and
\item  $a_i$ and $a_j$ are relatively prime whenever 
$D_i\cap D_j\neq \emptyset$, and
\end{enumerate}
\item a class   $B\in H^2(X,\z)$.
\end{enumerate}

There is a unique Seifert $\c^*$-bundle $f:Y\to X$
such that
\begin{enumerate}\setcounter{enumi}{3}
\item $f:Y\to X$ has  orbit invariants
$(a_i,b_i)$ along $D_i$, and
\item $f$ factors as 
$f:Y\stackrel{\pi}{\to}M_Y\stackrel{q}{\to}X$, 
where
\begin{enumerate}
\item $q:M_Y\to X$ is the unique $\c^*$-bundle with Chern class
 $aB+\sum b_i\frac{a}{a_i}[D_i]$ where $a=\lcm\{a_i\}$, and
\item $\pi:Y\to M_Y$ is an $a$-sheeted branched cover, branching  along
$q^{-1}(D_i)$ with multiplicity $a_i$.
\end{enumerate}
\end{enumerate}

Every Seifert $\c^*$-bundle $f:Y\to X$ has a unique such representation.

This representation defines the {\it Chern class} of a Seifert bundle
$$
c_1(Y/X):=B+\sum \tfrac{b_i}{a_i}[D_i]\in H^2(X,\q),
$$
where, for a divisor $D\subset X$, $[D]\in H^2(X,\z)$ denotes the
corresponding cohomology class. This notation extends to
$\q$-linear combinations of divisors by linearity.

If $X$ is projective and $H^2(X,\o_X)=0$ then 
every Seifert $\c^*$-bundle has a unique holomorphic 
Seifert $\c^*$-bundle structure. It satisfies $c_1(B_Y)=B$.
\end{say}

\section{The topology of 5-dimensional Seifert bundles}

\begin{notation}\label{top.seif.not}
 In this section, 
$X$ denotes a smooth, projective, simply connected
algebraic variety and
$D_1,\dots,D_n\subset X$ are smooth  divisors  intersecting transversally.
$f:L\to X$ denotes  a  Seifert bundle with orbit invariants
$(a_1,b_1),\dots, (a_n,b_n)$ along $D_1,\dots,D_n$.
 Set $a=\lcm(a_1,\dots,a_n)$.
\end{notation}

The main result, (\ref{seif.top.cor}) is for surfaces only, but
two of the intermediate steps hold in all dimensions.

\begin{prop}\label{top.of.Seif.H1} Notation as in (\ref{top.seif.not}).
 Assume that
\begin{enumerate}
\item the $[D_i]$ form part of a basis of $H_2(X,\z)$,
\item $a\cdot c_1(L)\in H^2(X,\z)$ is not divisible.
\end{enumerate}
Then $H_1(L,\z)=0$.

If, in addition,
\begin{enumerate}\setcounter{enumi}{2}
\item  $\pi_1(X\setminus(D_1\cup\cdots\cup D_n))$ is abelian,
\end{enumerate}
then $L$ is simply connected.
\end{prop}

\begin{prop}\label{top.of.Seif.H3}
Notation as in (\ref{top.seif.not}).
 Assume that  $\dim X=2$ and every $D_i$ is a rational curve.
Then $H^3(L,\z)$ is torsion free.
\end{prop}

\begin{prop}\label{top.of.Seif.w2} Notation as in (\ref{top.seif.not}).
 Assume that
\begin{enumerate}
\item the $a_i$ are odd, and
\item $w_2(X)\equiv a\cdot c_1(L)$ modulo 2.
\end{enumerate}
Then $w_2(L)$, the second
Stiefel--Whitney class of $L$ (cf.\ \cite[Sec.4]{milnor}), is zero.
\end{prop}

\begin{cor}\label{seif.top.cor}
 Notation as in (\ref{top.seif.not}). 
Assume that $\dim X=2$ and the conditions of (\ref{top.of.Seif.H1}),
(\ref{top.of.Seif.H3}) and (\ref{top.of.Seif.w2}) are all satisfied.
Then  $L$ is diffeomorphic to $(k-1)\#(S^2\times S^3)$
for $k=\dim H^2(X,\q)$.
\end{cor}

Proof. By \cite{smale},  a simply connected compact 5--manifold $L$
with vanishing second
Stiefel--Whitney class
 is uniquely determined by $H^3(L,\z)$.\qed
\medskip

The computation of $H_1(L,\z)$ relies on the following.

\begin{prop}\cite[4.6]{or-wa}\label{OW.h1.lem}
 Let $X$ be a complex manifold such that $H_1(X,\z)=0$ and let 
$D_1,\dots,D_n\subset X$ be smooth divisors intersecting transversally.
Let $f:L\to X$ be a  Seifert bundle with invariants
$(a_1,b_1,\dots, a_n,b_n,B)$. Then  $H_1(L,\z)$
is given by generators $k,g_1,\dots,g_n$ and relations
\begin{enumerate}
\item $a_ig_i+b_ik=0$ for $i=1,\dots,n$, and
\item $k(B\cap \eta)-\sum g_i([D_i]\cap \eta)=0$
for every $\eta\in H_2(X,\z)$.\qed
\end{enumerate}
\end{prop}

\begin{say}[Proof of (\ref{top.of.Seif.H1})]

It is enough to prove that the 
equations (\ref{OW.h1.lem}.1--2)
have only trivial solution modulo $p$ for every $p$.

If $p|a_i$ then $p\not\vert b_i$ so $a_ig_i+b_ik=0$ gives $k=0$. 
Since the $D_i$ form part of a basis of $H_2(X,\z)$, for every
$j$ there is an $\eta_j$ such that $[D_i]\cap \eta_j=\delta_{ij}$.
This implies that $g_j=0$ for every $j$.

If $p\not\vert a_i$ for every $i$ then $g_i=-(b_i/a_i)k$ makes sense
and the second equation, multiplied through by $a$, becomes
$$
k\cdot \bigl(aB+\sum b_i\tfrac{a}{a_i}[D_i]\bigr)\cap \eta=0
\qtq{for every $\eta\in H_2(X,\z)$.}
$$
By (\ref{top.of.Seif.H1}.2), $aB+\sum b_i\tfrac{a}{a_i}[D_i]=ac_1(L)$
is not zero modulo $p$, so for suitable $\eta$ we get $k=0$.

As in \cite[p.153]{or-wa}, the assumption (\ref{top.of.Seif.H1}.3)
 implies that 
$\pi_1(L)$ is abelian. Thus  $L$ is simply connected once $H_1(L,\z)=0$. \qed
\end{say}

\begin{say}[Proof of (\ref{top.of.Seif.H3})]

In order to compute the rest of the cohomology of $L$, 
we consider the Leray spectral sequence
$H^i(X,R^jf_*\z_L)\Rightarrow H^{i+j}(L,\z)$.
First
 we get some information about
the sheaf $R^1f_*\z_L$ and then about the groups $H^i(X,R^1f_*\z_L)$.

\begin{prop}\label{R^1.prop} Let $f:L\to X$ be a Seifert bundle.
\begin{enumerate}
\item There is a natural isomorphism
$\tau:R^1f_*\q_L\cong \q_X$.
\item There is a natural injection
$\tau:R^1f_*\z_L\into \z_X$ which is an isomorphism over the smooth locus.
\item If $U\subset X$ is connected then
$$
\tau(H^0(U,R^1f_*\z_L))= m(U)\cdot  H^0(U,\z)\cong m(U)\cdot \z,
$$
where $m(U)$ is the $\lcm$ of the multiplicities of
all fibers over $U$.
\end{enumerate}
\end{prop}

Proof. Pick $x\in X$ and a contractible neighborhood
$x\in V\subset X$. Then $f^{-1}(V)$ retracts to
$f^{-1}(x)\sim S^1$ and (together with the orientation)
this gives  a distinguished generator
$\rho\in H^1(f^{-1}(V), \z)$. This in turn determines
a  cohomology class $\frac1{m(x)}\rho\in H^1(f^{-1}(V), \q)$.
These normalized cohomology classes  are compatible with each other
and give a global section of $R^1f_*\q_L$. Thus
$R^1f_*\q_L=\q_X$ and we also obtain the
injection $\tau:R^1f_*\z_L\into \z_X$ as in (2). 

If $U\subset X$ is connected, a section of $b\in \z\cong H^0(U,\z_U)$ is in
$\tau(R^1f_*\z_L)$ iff $m(x)$ divides $b$ for every $x\in U$.
This is exactly (3).\qed

\begin{cor}\label{smooth.R^1.prop}
 Let  $f:L\to X$ 
be a Seifert bundle 
with orbit invariants $(a_i,b_i)$ along $D_i$.
Then there is an exact sequence
$$
0\to R^1f_*\z_L\stackrel{\tau}{\to} \z_X \to \sum_i \z_{D_i}/a_i\to 0.
$$
\end{cor}

Proof. Note that $a_i$ and $a_j$ are relatively prime if 
$D_i\cap D_j\neq \emptyset$. 
It is now clear that the kernel of 
$\z_X \to \sum_i \z_{D_i}/a_i$ has the same sections as
described in (\ref{R^1.prop}.3).\qed
\medskip

 The groups $H^i(X,R^1f_*\z_L)$
 sit in the long exact 
cohomology sequence of the short exact sequence of (\ref{smooth.R^1.prop}).
The  crucial  piece is
$$
 \sum_i H^1(D_i,\z/a_i)\to H^2(X,R^1f_*\z_L)\to H^2(X,\z).
$$
Thus $H^2(X,R^1f_*\z_L)$ is torsion free if $H^1(D_i,\z)=0$ for every $i$.

 Therefore, the $E_2$ term of the Leray spectral sequence
$H^i(X,R^jf_*\z_L)\Rightarrow H^{i+j}(L,\z)$  is
$$
\begin{array}{lcccc}
\z \quad & * & \quad\z^{k}\quad & * &\quad 
\z\\
\z & 0 & \z^{k} & 0  &\quad  \z.
\end{array}
$$
The spectral sequence degenerates at $E^3$ and we have only
two nontrivial differentials
$$
\delta_0: E^{0,1}_2\to E^{2,0}_2 \qtq{and} \delta_2:E^{2,1}_2\to E^{4,0}_2.
$$
In any case, $H^3(L,\z)\cong \ker \delta_2$
and so it is torsion free. \qed

Note also that if $H_1(L,\q)=0$ then $\delta_0$ is nonzero, hence
$\rank H^3(L,\q)=\rank H^2(L,\q)=\rank H^2(X,\q)-1$.
\end{say}

\begin{say}[Proof of (\ref{top.of.Seif.w2})]

Let $Y\supset L$ denote the corresponding Seifert $\c^*$-bundle.
$L$ is an orientable hypersurface, hence its normal bundle is trivial.
This implies that $w_i(L)=w_i(Y)|_L$. Since $w_2(Y)\equiv c_1(Y)\mod 2$
(cf.\ \cite[14-B]{milnor}),
it is enough to prove that $K_Y=-c_1(Y)$ is divisible by 2.

Let $E\to X$ be the unique holomorphic line bundle with
$c_1(E)=ac_1(L)$ 
with zero section $X\subset E$.
Let  $M:=E\setminus(\mbox{zero section})$ be
the corresponding $\c^*$-bundle.
By (\ref{seif.contr}.5.b), 
there is a branched covering  $\pi:Y\to M$ 
with branching multiplicities $a_i$. These are all odd,
 hence by
the Hurwitz formula, $K_Y\equiv \pi^* K_M\mod 2$.

By  the adjunction formula  
$K_X=K_E|_X+c_1(E)$. Thus, working in $H^2(X,\z/2)$, we get that 
$$
K_E|_X=K_X-c_1(E)=K_X-ac_1(L)\equiv w_2(X)-ac_1(L)\equiv 0 \qtq{modulo 2,}
$$
the last equality by (\ref{top.of.Seif.w2}.2).
The injection $X\into E$ is a homotopy equivalence,
thus $K_E$ and hence also $K_M$  are both divisible by 2.
\qed
\end{say}

\section{Einstein metrics on Seifert bundles}

 A  method  of Kobayashi \cite{kob} (see also \cite[9.76]{besse})
constructs Einstein metrics on 
circle bundles $M\to X$ from a K\"ahler--Einstein metric on $X$.
This was generalized to Seifert bundles $f:L\to X$ 
in  \cite{bg00}, but in this case one needs an 
orbifold  
K\"ahler--Einstein metric on $X$ and a Hermitian metric on $Y$.

\begin{defn} Let $f:Y\to X$ be a Seifert $\c^*$
bundle. A {\it Hermitian metric} on $Y$ is a
$C^{\infty}$ family of Hermitian metrics on the fibers.

We can also think of this as a degenerate
 Hermitian metric $h$ on the line bundle
$B_Y$. 
On  $X\setminus \cup D_i$ we get a Hermitian metric.
 On a local chart described in (\ref{stand.seif.defn}),
let $s(x_1,\dots,x_n)$ be 
a  generating section of $B_Y$.
Then we have the requirement
$$
h(s,s)=(\mbox{$C^{\infty}$-function})\cdot\prod (x_i\bar x_i)^{b_i/a_i}.
$$
The equivalence is clear since by  (\ref{stand.seif.defn})
we can think of $s$ as 
 $(\mbox{$C^{\infty}$-function})\cdot \prod x_i^{b_i/a_i}$.

This again leads to the Chern class equality
$c_1(B)+\sum \frac{b_i}{a_i}[D_i]=c_1(L)$.
\end{defn}

\begin{defn}  Let $X$ be a complex manifold and
$D_i\subset X$ smooth divisors intersecting transversally.
For every $i$ pick a natural number $a_i$.
Assume that $a_i$ and $a_j$ are relatively prime whenever 
$D_i\cap D_j\neq \emptyset$.

As in (\ref{stand.seif.defn}) and (\ref{seif.C*.defn}), 
for every $P\in X$ choose a neighborhood
$P\in U_P\subset X$ 
biholomorphic to  $D^n_{\mathbf x}$ which can be written as
a quotient 
$$
\pi_P:D^n_{\mathbf z}\to 
D^n_{\mathbf z}/\langle\phi_P\rangle\cong D^n_{\mathbf x}\cong U_P,
$$
where $\langle\phi_P\rangle$ is the cyclic group of order
$a_P=\prod_{P\in D_i} a_i$
and $\pi_P$ branches exactly along
the divisors 
 $D_i\cap U_P$  with  multiplicity $a_i$. 

Thus 
the action is 
$$
\phi_P: (z_1,\dots,z_n)\mapsto 
(\epsilon_1z_1,\dots,\epsilon_nz_n),
$$
where 
 $\epsilon_j$ is a primitive $a_{i_j}$th root of unity
for $D_{i_j}\cap U_P\neq \emptyset$
and we set $\epsilon_j=1$ for those coodinates that do
not correspond to any $D_i$.

An {\it orbifold  
Hermitian metric} on $(X,\sum (1-\frac1{a_i})D_i)$
is a Hermitian metric $h$ on $X\setminus \cup D_i$
such that $\pi_P^*h$ extends to a 
Hermitian metric on $D^n_{\mathbf z}$ for every $P\in X$.

 On a local chart described in (\ref{stand.seif.defn})
this means a metric
$\sum h_{ij}dx_i\otimes d\bar x_j$ 
defined on $D^n_{\mathbf x}\setminus (\prod x_i=0)$
such that
$\sum a_ia_jh_{ij}z_i^{(a_i-1)}\bar z_j^{(a_j-1)}dz_i\otimes d\bar z_j$
is a metric on  $D^n_{\mathbf z}$.
Thus the  orbifold canonical class is $K_X+\sum (1-\frac1{a_i})D_i$.

An {\it orbifold  
K\"ahler--Einstein metric} is now defined as usual.
\end{defn}

\begin{thm} \cite{kob, bg00}\label{exst.of.E.metric}
Let $f:L\to X$ be a
Seifert bundle with orbit invariants $(a_1,b_1),\dots, (a_n,b_n)$.
 $L$ admits an $S^1$-invariant Einstein metric 
with positive constant if the following hold.
\begin{enumerate}
\item The orbifold canonical class
$K_X+\sum (1-\frac1{a_i})D_i$ is anti ample and 
there is a K\"ahler--Einstein metric on $(X,\sum (1-\frac1{a_i})D_i)$.
\item The Chern class of  $L$ is a negative multiple of
$K_X+\sum (1-\frac1{a_i})D_i$. \qed
\end{enumerate}
\end{thm}

 The main impediment to apply
(\ref{exst.of.E.metric}) is the current shortage of
existence results for K\"ahler--Einstein metrics on orbifolds.
We use the following sufficient algebro--geometric condition.
There is every reason to expect that it is very far from being optimal,
but it does provide a large selection of  examples.

In this paper we use (\ref{nadel.thm}) only for surfaces.
The concept {\it klt} is defined in (\ref{klt.etc.defn}).

\begin{thm}\cite{nadel, dk}
\label{nadel.thm}  Let $X$ be an $n$-dimensional 
 compact complex manifold and
$D_i\subset X$ smooth divisors intersecting transversally.

 Assume that $-(K_X+\sum (1-\frac1{a_i})D_i)$ is ample and 
there is an $\epsilon>0$ such that 
$$
(X,\tfrac{n+\epsilon}{n+1}F+\sum (1-\tfrac1{a_i})D_i) \qtq{is klt}
\eqno{(\ref{nadel.thm}.1)}
$$
for every positive $\q$-linear combination of divisors $F=\sum f_iF_i$
such that $[F]= -[K_X+\sum (1-\frac1{a_i})D_i]\in H^2(X,\q)$.

 Then there is  orbifold  K\"ahler--Einstein metric
on $(X,\sum (1-\frac1{a_i})D_i)$. 

If $\sum a_iD_i$ is invariant under a compact group $G$ of biholomorphisms
of $X$, then it is sufficient to check (\ref{nadel.thm}.1) 
for $G$-equivariant divisors $F$.
\qed
\end{thm}

\begin{defn}(cf. \cite[2.34]{kmbook}\label{klt.etc.defn}
 Let $X$ be a complex manifold  and $D$ an effective  $\q$-divisor on $X$.
Let $g:Y\to X$ be any proper bimeromorphic  morphism, $Y$ smooth.
 Then there is a unique
$\q$-divisor $D_Y=\sum e_iE_i$ on $Y$ such that 
$$
K_Y+D_Y\equiv g^*(K_X+D)\qtq{and}  g_*D_Y=D.
$$
We say that $(X,D)$ is {\it klt}
(short for Kawamata log terminal)
if $e_i< 1$ for every $g$
and for every $i$.
\end{defn}

It is quite hard to check using the above 
 definition if a pair $(X,D)$ is klt or not.
 For surfaces, there are reasonably sharp
multiplicity conditions  which ensure that a given pair $(X,D)$ 
is klt. These conditions are not  necessary, but
they seem to apply in most cases of interest to us.

\begin{lem}(cf.\ \cite[4.5 and 5.50]{kmbook})
\label{how.to.check}
Let $S$ be a smooth surface and $D$ an effective $\q$-divisor.
Then   $(S,D)$ is klt if the following 
conditions  are satisfied.
\begin{enumerate}
\item   $D$ does not contain an
irreducible component with coefficient  $\geq 1$.
\item  For every  point $P\in S$, either
\begin{enumerate}
\item $\mult_PD\leq 1$, or
\item we can write $D=cC+D'$ where $C$ is a curve through $P$, smooth at $P$,
$D'$ is effective not containing $C$,
and the local intersection number $(C\cdot_P D')<1$.\qed
\end{enumerate}
\end{enumerate}
\end{lem}

\section{Log Del Pezzo surfaces with large $H^2$}

In this section we construct smooth log Del Pezzo surfaces with 
a K\"ahler--Einstein  metric and $\rank H^2$ arbitrarily large.
Methods to construct log Del Pezzo surfaces with large
$H^2$ are given in \cite{keel-mck, mck}.
For ordinary smooth Del Pezzo surfaces the rank of $H^2$
is at most 9.

\begin{exmp} Start with $\p^1\times \p^1$ with coordinate
projections $\pi_1,\pi_2$. Let $C_1\subset \p^1\times \p^1$
be a fiber of $\pi_2$ and  $C_2\subset \p^1\times \p^1$
the graph of a degree 2 morphism $\p^1\to \p^1$.
Thus $\pi_1:C_2\to \p^1$ has degree 1 and 
 $\pi_2:C_2\to \p^1$ has degree 2. Pick $k$ points
$P_1,\dots,P_k\in C_2\setminus C_1$
and blow them up to obtain a surface $h:S_k\to \p^1\times \p^1$.
Let $C_i'\subset S_k$ denote the birational transform of $C_i$.
Note that $(C_1')^2 =0$, $(C_1'\cdot C_2')=2$ and  $(C_2')^2 =4-k$.

$C_1+C_2\in |-K_{\p^1\times \p^1}|$, thus
$C_1'+C_2'\in |-K_{S_k}|$. 
\end{exmp}

\begin{lem}\label {ldp.amp.lem}
Let $a_i$ be rational numbers with   $a_1>\frac{k-4}{2}a_2>0$. Then
a large multiple of $a_1C_1'+a_2C_2'$  determines a birational morphism
$S_k\to \bar S_k$. The positive dimensional fibers are exactly the
birational transforms of those fibers of $\pi_2$ which intersect
$C_2$ in two points of the set  $P_1,\dots,P_k$.
\end{lem}

Proof. The conditions ensure that $(a_1C_1'+a_2C_2')\cdot C_2'>0$,
thus $a_1C_1'+a_2C_2'$ is nef and big.

For $c>a_1,a_2$ we can write
$$
a_1C_1'+a_2C_2'\equiv -c\left( K_{S_k}+(1-a_1/c)C_1'+(1-a_2/c)C_2'\right).
$$
The Base point free
theorem (cf.\ \cite[3.3]{kmbook}) applies and so
a large multiple of $a_1C_1'+a_2C_2'$  determines a birational morphism.

 The positive dimensional fibers are exactly those curves which have
zero intersection number with $a_1C_1'+a_2C_2'$.
The projection of such a curve to $\p^1\times \p^1$
is thus a curve $B$ which intersects $C_1+C_2$ only at the points
$P_1,\dots,P_k$. Since $B$ is disjoint from $C_1$,
it is the union of fibers of $\pi_2$.\qed

\begin{lem}\label{ldp.top.lem} Notation as above. Assume that $k\geq 5$
and choose natural numbers $m_1,m_2\geq 2$ satisfying
$m_2>\frac{k-4}{2} m_1$. Assume that no two of the $P_i$ are on 
the same fiber of $\pi_2$. Then
\begin{enumerate}
\item $S_k$ is smooth, its  Picard number  is $k+2$.
\item $C_1'$ and $C_2'$ are smooth rational curves intersecting transversally
and they form part of a basis of $H^2(S^*_k,\z)$.
\item $-(K_{S_k}+(1-\tfrac1{m_1})C_1'+(1-\tfrac1{m_2})C_2')$ is ample.
\item  $\pi_1(S_k\setminus(C_1'+C_2'))=1$.
\end{enumerate}
\end{lem}

Proof. The canonical class of $S_k$ is $-C_1'-C_2'$, thus
$$
-(K_{S_k}+(1-\tfrac1{m_1})C_1'+(1-\tfrac1{m_2})C_2')=
\tfrac1{m_1}C_1'+\tfrac1{m_2}C_2',
$$
and this is ample by (\ref{ldp.amp.lem}).

By looking at the first projection
$\pi_1$ we see that
$\p^1\times \p^1\setminus(C_1+C_2)$ contains an open subset
$\c^*\times \c^*$. 
Thus  $\pi_1(S_k\setminus(C_1'+C_2'))$ is abelian
and is generated by small loops around the $C_i'$. 

Any of the exceptional curves of $S_k\to \p^1\times \p^1$
intersects $C_2'$ transversally in one point and is disjoint from
$C_1'$. Similarly, the birational transform
of $\p^1\times \{P_i\}$ (for any $i$) intersects $C_1'$ transversally in 
one point and is disjoint from
$C_2'$. These show that $C_1'$ and $C_2'$ 
 form part of a basis of $H^2(S_k,\z)$ and that 
the small loops around the $C_i'$
are contractible in $S_k\setminus(C_1'+C_2')$, thus
 $\pi_1(S_k\setminus(C_1'+C_2'))=1$.
\qed
\medskip

The next result gives information about klt divisors
on $S_k$.

\begin{lem}\label{ldp.klt.1pt}  Notation as above and assume that $k\geq 5$.
Let $D\subset S_k$ be an effective $\q$-divisor
such that $[D]=[b_1C_1'+b_2C_2']$ for some $0\leq b_i<1/2$.
Then $(S_k,\frac12C_1'+\frac12C_2'+D)$ is klt, except possibly at 
one, but not both,  of the
two intersection points $C_1'\cap C_2'$.
\end{lem}

Proof. Assume that $\frac12C_1'+\frac12C_2'+D$ is not klt at a
point $Q\in S_k$. 
If $Q$ is not on any of the exceptional curves of $h$
then 
$h_*(\frac12C_1'+\frac12C_2'+D)$ is also 
not klt at $h(Q)\in \p^1\times \p^1$. By assumption
$[h_*D]$ is cohomologous to the sum of the two lines on $\p^1\times \p^1$
with coefficients less than 1;  these are always klt.
Thus
$\frac12C_1'+\frac12C_2'+D$ is klt outside 
$h^{-1}(C_1+C_2)$.

Write 
$$
\tfrac12C_1'+\tfrac12C_2'+D=(\tfrac12+d_1)C_1'+(\tfrac12+d_2)C_2'+D',
$$
where $D'$ does not contain the $C_i'$. 
Then $D'\equiv (b_1-d_1)C_1'+(b_2-d_2)C_2'$ and
$d_i\leq b_i$ since
 otherwise we would obtain that
$$
[D']=\pm[\alpha C_1'-\beta C_2']\qtq{with $\alpha,\beta >0$.}
$$
Both are impossible as they lead to  a negative intersection number with
one of the $C_i'$. In particular, $\frac12+d_i<1$, so
the $C_i'$ appear in $\frac12C_1'+\frac12C_2'+D$
with coefficient less than 1. Thus
$\frac12C_1'+\frac12C_2'+D$ satisfies the condition (\ref{how.to.check}.1)
for $C_1'$ and $C_2'$.

In order to  check  (\ref{how.to.check}.2.b) along $C_1'$
(resp.\  $C_2'$), we study  the restriction
$(\frac12+d_2)C_2'+D'|_{C_1'}$ (resp. $(\frac12+d_1)C_1'+D'|_{C_2'}$)
is klt.
We compute the  intersection numbers 
$$
\begin{array}{rcl}
\deg D'|_{C_1'}&=&(D'\cdot C_1')
=2(b_2-d_2)<1,\qtq{and}\\
\deg D'|_{C_2'}&=&(D'\cdot C_2')
=2(b_1-d_1)-(k-4)(b_2-d_2)<1.
\end{array}
$$
In both cases, the remainder of the restrictions consists of the
2 intersection points $Q_1+Q_2=C_1'\cap C_2'$, each with coefficient
$\frac12+d_2$ (resp. $\frac12+d_1$). 
Therefore
$\deg \left((\frac12+d_{3-i})C_{3-i}'+D'|_{C_i'}\right)<2$
and the restriction contains both $Q_1$ and $Q_2$ with coefficient
bigger than $\frac12$.

Thus  every point occurs
with coefficient less than 1, except possibly for $Q_1$ and $Q_2$.
  Moreover, only one of these two points can have
coefficient at least 1. 

We are left to understand what happens along the exceptional curves
of $h$. Let $E$ be such a curve and write
$D'=eE+D''$ where $D''$ does not contain $E$. Then
$$
\begin{array}{l}
e\leq (D'\cdot C'_2)\leq 2(b_1-d_1)-(k-4)(b_2-d_2)<1,
\qtq{and so}\\
(D''\cdot E)=  (D'\cdot E)+e\leq 2(b_1-d_1)-(k-5)(b_2-d_2)<1.
\end{array}
$$
The first inequality shows (\ref{how.to.check}.1) for $E$
and the second shows that (\ref{how.to.check}.2.b) holds at every point
of $E$.
\qed
\medskip

\begin{rem} In the above proof we could  have used the
Connectedness theorem
(cf.\ \cite[5.48]{kmbook}), which implies that
the set of points where 
$(S_k,\frac12C_1'+\frac12C_2'+D)$ is not klt is connected.
\end{rem}

This is, however, not enough to obtain a K\"ahler--Einstein
  metric on $S_k$.
To achieve this, we make a special choice of the points $P_i$.
It is easiest to write down everything by equations.

\begin{defn}\label{symmetric.constr.defn}
Choose homogeneous coordinates
$\p^1_{(s:t)}\times \p^1_{(u:v)}$ and pick
$C_1=(u=v)$ and $C_2=(s^2u=t^2v)$. The two intersection points $Q_i$ are
$(\pm1:1,1:1)$. The involution
$\tau: (s:t,u:v)\mapsto (-t:s,v:u)$ fixes that $C_i$
and interchanges the two points $Q_i$.

If $k=2m$ is even, pick $0<c_1<\cdots <c_m<1$
and for the $P_i$  choose the $2m$ points
$$
(c_i:1, c_i^2:1)\qtq{and} (-1:c_i,1:c_i^2),
$$
to obtain a surface $S^*_k$ of Picard number $k+2$.

If $k=2m+1$ is odd, pick $0<c_1<\cdots <c_m<1$
and  for the $P_i$ choose the $2m+1$ points
$$
(c_i:1, c_i^2:1),(-1:c_i,1:c_i^2)\qtq{and} (\sqrt{-1}:1,-1,1),
$$
to obtain a surface $S^*_k$ of Picard number $k+2$.

In both cases, the involution $\tau$ lifts to an involution
on $S^*_k$, 
again denoted by $\tau$.
\end{defn}

\begin{lem}\label{ldp.klt.symm}
 Notation as above and assume that $k\geq 5$.
Let $D\subset S^*_k$ be an effective $\tau$-invariant $\q$-divisor
such that  $[D]=[b_1C_1'+b_2C_2']$ for some $0\leq b_i<1/2$.
Then $(S^*_k, \frac12C_1'+\frac12C_2'+D)$ is klt.
\end{lem}

Proof. We already know that  $\frac12C_1'+\frac12C_2'+D$ is klt,
except possibly at the
two intersection points $C_1'\cap C_2'$. These are interchanged by
$\tau$, so if  $\frac12C_1'+\frac12C_2'+D$ is not klt, then
it is not klt at exactly these 2 points. We have seen
in (\ref{ldp.klt.1pt})
 that it can not fail to be klt at both of these points. \qed

\begin{cor}\label{ldp.ke.cor} Notation as above. Assume that $k\geq 5$
and choose natural numbers $m_1,m_2\geq 2$ satisfying
$m_2>\frac{k-4}{2} m_1$. Then 
 $(S^*_k, (1-\tfrac1{m_1})C_1'+(1-\tfrac1{m_2})C_2')$ has an
orbifold K\"ahler--Einstein  metric.
\end{cor}

Proof. Set $\Delta=(1-\tfrac1{m_1})C_1'+(1-\tfrac1{m_2})C_2'$.
Then $-(K_{S^*_k}+\Delta)=\tfrac1{m_1}C_1'+\tfrac1{m_2}C_2'$
is ample by (\ref{ldp.amp.lem}).  For the existence of 
a K\"ahler--Einstein  metric,
we use the criterion (\ref{nadel.thm}) with $\frac23+\epsilon=\frac34$. 
Let $D$ be any effective $\tau$-invariant divisor
numerically equivalent to $-(K_{S^*_k}+\Delta)$.
Then
$$
\Delta+\tfrac34 D\equiv \tfrac12C_1'+\tfrac12C_2'+
(\tfrac12-\tfrac1{m_1})C_1'+(\tfrac12-\tfrac1{m_2})C_2'+\tfrac34D.
$$
Note that
$$
(\tfrac12-\tfrac1{m_1})C_1'+(\tfrac12-\tfrac1{m_2})C_2'+\tfrac34D
\equiv 
(\tfrac12-\tfrac1{4m_1})C_1'+(\tfrac12-\tfrac1{4m_2})C_2'.
$$
The assumptions of (\ref{ldp.klt.symm}) are satisfied
and so $(S^*_k, \Delta+\tfrac34 D)$ is klt. Thus
$(S,\Delta)$ has an orbifold K\"ahler--Einstein  metric.
\qed

\section{Einstein metrics on $k\#(S^2\times S^3)$}

Let $C_1,C_2\subset \p^1\times \p^1$ be as in
(\ref{symmetric.constr.defn}).
For $k\geq 6$ let $M_{k-1}$ be the moduli space of $k-1$ distinct points
in $C_2\setminus C_1$. Its dimension is $k-1$.
We can also think of $M_{k-1}$ as
parametrizing surfaces of type $S_{k-1}$ obtained by blowing up
these points.

Fix  $k\geq 6$
and relatively prime odd numbers $m_1,m_2>2$ satisfying 
$m_2>\frac{k-5}{2}m_1$.
We can then further view  $M_{k-1}$ as
parametrizing orbifolds 
$$
(S_{k-1},(1-\tfrac1{m_1})C_1'+(1-\tfrac1{m_2})C_2').
$$

Let $f:L\to S_{k-1}$ be the Seifert bundle with
orbit invariants $(m_1,1)$ along $C_1'$,
$(m_2,1)$ along $C_2'$ and with the trivial line bundle as
$B_L$.

The number $a$ in (\ref{seif.contr}.6) is $m_1m_2$.
The Chern class of $L$ is
$$
c_1(L)=\tfrac1{m_1}[C_1']+\tfrac1{m_2}[C_2']
=-\left(K_{S_{k-1}}+(1-\tfrac1{m_1})C_1'+(1-\tfrac1{m_2})C_2')\right).
$$
Furthermore, $ ac_1(L)=m_2[C_1']+m_1[C_2']$
 is not divisible since $m_1,m_2$ are 
relatively prime and computing modulo 2 
$$
ac_1(L)\equiv -[C_2']-[C_1']=[K_{S_{k-1}}]\equiv w_2(S)
$$
 since the $m_i$ are odd.

Using  (\ref{ldp.top.lem}) as well, we see that the conditions of
(\ref{seif.top.cor}) are all satisfied, thus 
any such $L$ is 
diffeomorphic to $k\#(S^2\times S^3)$.

The surfaces of type $S^*_{k-1}$ considered in (\ref{ldp.ke.cor})
form a subset of $M_{k-1}$, and for these surfaces
we have proved the existence of an orbifold K\"ahler--Einstein  metric.
The existence of an orbifold K\"ahler--Einstein  metric
is an open condition (in the Euclidean topology), thus we obtain:

\begin{claim} For $k\geq 6$ and odd $m_1,m_2>2$ satisfying 
$m_2>\frac{k-5}{2}m_1$ there is an open subset
$U(k-1,m_1,m_2)\subset M_{k-1}$ such that  all orbifolds
$(S_{k-1},(1-\tfrac1{m_1})C_1'+(1-\tfrac1{m_2})C_2')$
corresponding to a point in $U(k-1,m_1,m_2)$
 have an
orbifold K\"ahler--Einstein  metric.
\end{claim}

By (\ref{exst.of.E.metric}), for any surface corresponding to a point in
$U(k-1,m_1,m_2)$ we obtain an Einstein metric on $L$. The space 
$U(k-1,m_1,m_2)$ has complex dimension $k-1$ hence real dimension $2k-2$.
\qed

\begin{complem} As explained in \cite{bg00} (see also \cite{bgk}),
the Einstein metrics constructed this way have additional good properties.
\begin{enumerate}
\item The connected component of the isometry group of the  metric
is $S^1$.
\item  All these metrics are
{\it Sasakian-Einstein}.
\item Two metrics constructed from the data
$(m_1^i,m_2^i, P_1^i,\dots,P_{k-1}^i)$  for $i=1,2$ are isometric
iff  $m_1^1=m_1^2, m_2^1=m_2^2$ and the point set
$\{P_1^1,\dots,P_{k-1}^1\}$ can be mapped to either
$\{P_1^2,\dots,P_{k-1}^2\}$ or to its conjugate
by an automorphism of $\p^1\times \p^1$ fixing $C_1$ and $C_2$.
(Such  automorphisms form a group of order 4.)
\end{enumerate}
\end{complem}

\begin{ack} I thank Ch.\ Boyer, K.\ Galicki and J.\ McKernan
for many useful conversations and e-mails.
Research
was partially supported by the NSF under grant number
DMS-0200883. 
\end{ack}

\bibliography{refs}

\end{document}